%% file: main.tex
\documentclass[12pt,a4paper,leqno]{article}

\usepackage{latexsym}
\usepackage{amssymb,exscale}
\usepackage[centertags]{amsmath}
\usepackage{amsthm}
\usepackage{graphicx}

\input{preamble}

\begin{document}

\title{Random complex zeroes, \rm{III}.
\\ Decay of the hole probability}
\author{Mikhail Sodin and Boris Tsirelson}
\date{}

\maketitle

\begin{abstract}
The `hole probability' that a random entire function
\[
\psi(z) = \sum_{k=0}^\infty \zeta_k \frac{z^k}{\sqrt{k!}} \, ,
\]
where $\zeta_0, \zeta_1, \dots $ are Gaussian i.i.d.\ random
variables, has no zeroes in the disc of radius $r$ decays as $ \exp
(-c r^4)$ for large $r$.
\end{abstract}

We consider the (random) set of zeroes of a random entire function
$\psi_\omega:~\C\to \C$,
\begin{equation}\label{eq0.1}
\psi (z, \omega) =
\sum_{k=0}^\infty \zeta_k(\omega) \frac{z^k}{\sqrt{k!}}\,,
\end{equation}
where $\zeta_k$, $k=0, 1, 2$, ... are independent standard complex-valued
Gaussian random variables, that is the distribution $\mathcal
N_{\C}(0,1)$ of each $\zeta_k$ has the density $\pi^{-1} \exp(-|w|^2)$
with respect to the Lebesgue measure $m$ on $\C$. This model is
distinguished by invariance of the distribution of zero points with
respect to the motions of the complex plane
\[
z\mapsto az+b, \qquad |a|=1, \ b\in\C\,,
\]
see \cite{SoTsi} for details and references.

Given large positive $r$, we are interested here in the `hole probability'
that $\psi$ has no zeroes in the disc of radius $r$
\[
p(r) = \mathbb P\big(\psi (z,\cdot)\ne 0, \ |z|\le r \big)\,.
\]
It is not difficult to show that $ p(r) \le \exp(-\const\, r^2)$,
see the Offord-type estimate in \cite{So}. Yuval Peres told one of
us that the recent work \cite{PV} led to conjecture that the
actual hole probability might have a faster decay. In this note,
we confirm this conjecture and prove
\begin{theorem}\label{thm1}
$\exp(-C r^4)
\le p(r) \le \exp(-c r^4)$.
\end{theorem}
Throughout, by $c$ and $C$ we denote various positive numerical constants
whose values can be different at each occurrence.

It would be interesting to check whether there exists the limit
\[
\lim_{r\to\infty} \frac{\log^- p(r)}{r^4}\,,
\]
and (if it does) to find its value.

The lower bound in Theorem~\ref{thm1} will be obtained in
Section~1 by a straightforward construction. The upper bound in
Theorem~\ref{thm1} follows from a large deviation estimate which
has an independent interest.
\begin{theorem}\label{thm2}
Let $n(r)$ be a number of random zeroes in the disc $\{|z|\le r\}$. Then
for any $\de\in (0, \frac14]$ and $r\ge 1$
\begin{equation}\label{eq0.3}
\mathbb P \bigg( \bigg| \frac{n(r)}{r^2} - 1 \bigg| \ge \de \bigg)
\le \exp(-c(\de)r^4)\,.
\end{equation}
\end{theorem}
Throughout, by $c(\de)$ we denote various positive constants which
depend on $\de$ only. Since our argument is too crude to find a
sharp constant $c(\de)$ in \eqref{eq0.3}, we freely change the
values of $c(\de)$ from line to line.

There is a fruitful analogy between random zero sets and one
component Coulomb system which consists of charged particles of one
sign in $\R^2$ embedded in a uniform background of the opposite
sign (see \cite{FH} and references therein). Theorems~\ref{thm1}
and \ref{thm2} are consistent with the corresponding results for
Coulomb systems \cite{JLM}.

\medskip\par\noindent{\sc Acknowledgment. }
Yuval Peres brought our attention to the problem considered here.
F\"edor Nazarov spotted an error in the first draft and suggested
how to fix it. We thank both of them.

\section{Proof of the lower bound in Theorem~\ref{thm1}}

In what follows, we frequently use two elementary facts:
if $\zeta$ is a standard complex Gaussian variable, then
\begin{equation}\label{fact1}
\mathbb P(|\zeta|\ge \lambda) = \frac1{\pi} \iint_{|w|\ge \lambda}
e^{-|w|^2}\, dm(w) = \int_{\lambda^2}^\infty e^{-t}\, dt =
e^{-\lambda^2},
\end{equation}
and for $\lambda\le 1$
\begin{equation}\label{fact2}
\mathbb P(|\zeta|\le \lambda) = 1-e^{-\lambda^2} = \lambda^2 -
\frac{\lambda^4}{2!} + \, ...\, \in \bigg[ \frac{\lambda^2}2,
\lambda^2 \bigg].
\end{equation}

By $\Omega_r$ we denote the following event: (i) $ |\zeta_0|\ge
2$; (ii) $|\zeta_k|\le \exp (-2r^2)$ for $1\le k \le 48r^2$; and
(iii) $|\zeta_k|\le 2^k$ for $k>48r^2$. Since $\zeta_k$ are
independent,
\[
\mathbb P(\Omega_r) = \mathbb P\hbox{(i)}\cdot \mathbb
P\hbox{(ii)}\cdot \mathbb P\hbox{(iii)}\,.
\]
Evidently, the first and third factors on the RHS are $\ge
\const$. By (\ref{fact2}), the probability of the event
$|\zeta_k|\le \exp(-2r^2)$ is $\ge \frac12 \exp(-4r^2)$. Since the
events within the second group are independent, the probability of
all of them to happen is $ \ge \big( \frac12 \exp(-4r^2) \big)^{48r^2}
= \exp(-192r^4-Cr^2)$. Thus, $\mathbb P(\Omega_r) \ge \exp(-Cr^4)$.

Now, we show that for $\omega\in\Omega_r$ the function $\psi $
does not vanish in the disc $\{|z|\le r\}$. For such $z$ and $\omega$
we have
\[
|\psi (z)| \ge |\zeta_0| - \sum_{1\le k\le 48 r^2} |\zeta_k|
\frac{r^k}{\sqrt{k!}} - \sum_{k> 48r^2} |\zeta_k|
\frac{r^k}{\sqrt{k!}} = |\zeta_0| - \sum\nolimits' -
\sum\nolimits''.
\]
Then
\begin{eqnarray*}
\sum\nolimits' &\stackrel{\hbox{(ii)}}\le& e^{-2r^2} \sum_{1\le k
\le 48r^2} \frac{r^k}{\sqrt{k!}} \\
&\le& e^{-2r^2} \, \sqrt{48r^2}\cdot\sqrt{\sum_{1\le k\le 48r^2}
\frac{r^{2k}}{k!}}  < 7r\, e^{-2r^2 + 0.5r^2} < e^{-r^2} <
\frac12\,,
\end{eqnarray*}
if $r$ is sufficiently large. At the same time,
\[
\sum\nolimits'' \stackrel{\hbox{(iii)}}\le \sum_{k>48r^2}
\frac{2^k}{\sqrt{k!}} \bigg( \frac{k}{48} \bigg)^{k/2}
< \sum_{k>48r^2} \bigg(
\frac{k}{12} \cdot \frac{e}{k} \bigg)^{k/2} < \sum_{k\ge 1}
2^{-k} = \frac12
\]
(we used inequality $k!> \big( \frac{k}{e} \big)^k$ which
follows from Stirling's formula). Putting both estimates together,
we get
\[
|\psi (z)| \ge |\zeta_0| - 1 \stackrel{\hbox{(i)}}\ge 1\,, \qquad
|z|\le r\,,
\]
proving that $\psi $ does not vanish in the closed disc $\{|z|\le
r\}$ for $\omega\in\Omega_r$.

\section{Large deviations of $\log M(r, \psi) - r^2/2$}

Let $\psi$ be the random entire function \eqref{eq0.1} and let
$M(r,\psi)=\max_{|z|\le r}|\psi(z)|$. In this section we shall
prove the following
\begin{lemma}\label{lem2.1}
Given $\de\in (0, \frac14]$ and $r\ge 1$,
\[
\mathbb P \bigg( \bigg| \frac{\log M(r, \psi)}{r^2}
- \frac12 \bigg|\ge \de \bigg) \le \exp(-c(\de) r^4) \, .
\]
\end{lemma}

The proof is naturally split into two parts.
First we show that
\begin{equation}\label{eq2.2}
\mathbb P \bigg( \frac{\log M(r, \psi)}{r^2} \ge \frac12 + \de
\bigg) \le \exp(-c(\de) r^4)\,,
\end{equation}
and then that
\begin{equation}\label{eq2.3}
\mathbb P \bigg( \frac{\log M(r, \psi)}{r^2} \le \frac12 - \de
\bigg) \le \exp(-c(\de) r^4)\,.
\end{equation}

\begin{proof}[Proof of \eqref{eq2.2}]
We use an argument similar to the one used in Section~1. We have
\[
M(r, \psi) \le \bigg( \sum_{0\le k<4er^2} + \sum_{k\ge 4er^2}
\bigg) |\zeta_k| \frac{r^k}{\sqrt{k!}} = \sum\nolimits_1 +
\sum\nolimits_2\,.
\]
Consider the event $A_r$ which consists of such $\omega$'s that
(i) $|\zeta_k|\le \exp(2\de r^2/3)$ for
$0\le k<4er^2$; (ii) $|\zeta_k|\le (\sqrt 2)^k$ for $k\ge 4er^2$.
If $A_r$ occurs and $r$ is sufficiently large, then
\begin{eqnarray*}
\sum\nolimits_1^2 &\le& \bigg( \sum_{0\le k <4er^2}
|\zeta_k|^2 \bigg) \cdot \bigg( \sum_{0\le k < 4er^2}
\frac{r^{2k}}{k!}\bigg) \\
&\stackrel{\hbox{(i)}}\le& 4er^2\cdot \exp(4\de r^2/3 + r^2) <
\exp\bigg( \Big(1 + \frac53 \de\Big)r^2\bigg)\,,
\end{eqnarray*}
and
\[
\sum\nolimits_2 \stackrel{\hbox{(ii)}}\le \sum_{k\ge 4er^2}
|\zeta_k| \bigg( \frac{k}{4e} \cdot \frac{e}{k}\bigg)^{k/2} \le
\sum_{k\ge 4er^2} \frac{(\sqrt 2)^k}{2^k} \le 1.
\]
Thus
\[
M(r, \psi) \le \exp\bigg( \Big(\frac12 + \de\Big)r^2 \bigg)\,.
\]

It remains to estimate the probability of the complementary set
$A_r^c= \Omega\setminus A_r$. If $A_r^c$ occurs, then at least one of
the following happens:
$\exists k\in [0, 4er^2)$: $|\zeta_k|\ge \exp(\frac23 \de r^2)$, or
$\exists k\in [4er^2, \infty)$: $|\zeta_k|\ge (\sqrt 2)^k$. Therefore
\[
\mathbb P(A_r^c) \le 4er^2 \exp\bigg( \! -\exp\Big(\frac43 \de
r^2\Big) \bigg) + \sum_{k\ge 4er^2} \exp\big(-2^k \big) <
\exp\big( -\exp(\de r^2)\big)
\]
provided that $r\ge r_0(\de)$.  This is much stronger than
\eqref{eq2.2}.
\end{proof}

\begin{proof}[Proof of \eqref{eq2.3}]
Suppose that
\begin{equation}\label{eq2.4}
\log M(r, \psi) \le \bigg(\frac12 - \de\bigg)r^2\,.
\end{equation}
Then we use Cauchy's inequalities and Stirling's formula:
\begin{eqnarray*}
|\zeta_k| &=& \frac{|\psi^{(k)}(0)|}{\sqrt{k!}} \le \sqrt{k!}\,
\frac{M(r, \psi)}{r^k} \\
&\le& Ck^{1/4}\, \exp\bigg( \frac{k}{2}\log k - \frac{k}2 +
\Big(\frac12-\de \Big)r^2 - k\log r \bigg)\,.
\end{eqnarray*}
Observe that the exponent equals
\[
\frac{k}2 \bigg( (1-2\de) \frac{r^2}{k} - \log \frac{r^2}{k} -
1\bigg) \, .
\]
We note that $ (1-2\de) \frac{r^2}{k} - \log \frac{r^2}{k} - 1 < -\de
$ when $ r^2/k $ is close enough to $ 1 $.
Whence, for $(1-\epsilon)r^2 \le k \le r^2$,
\[
|\zeta_k| \le Ck^{1/4} \, \exp\bigg( -\frac{k\de}2\bigg)\,.
\]
By (\ref{fact2}), the probability of this event is $\le \exp
\big(- c(\de)k \big)$. Since $\zeta_k$ are independent,
multiplying these probabilities, we see that
\[
\exp\bigg( - c(\de) \sum_{(1-\epsilon)r^2 \le k \le r^2} k \bigg)
= \exp\big( -c_1(\de) r^4 \big)
\]
is an upper bound for the probability that event \eqref{eq2.4} occurs.
\end{proof}

\section{Mean lower bound for $\log|\psi(z)| - |z|^2/2$}

Lemma~\ref{lem2.1} gives us a sharp upper bound for the `random potential'
$ \log|\psi(z)| - \frac12|z|^2 $ when $\omega$ does not belong to an exceptional
set in the probability space. Here, we give a mean lower bound for this
potential.

\begin{lemma}\label{lem3.1}
Given $\de\in (0, \frac14]$ and $r\ge1$,
\[
\mathbb P\bigg( \frac1{r^2} \int_{r\T} \log|\psi|\, d\mu \le \frac12 -\de
\bigg) \le \exp(-c(\de)r^4)\,.
\]
\end{lemma}

Here, we denote by $r\T$ the circle $\{|z|=r\}$, $\mu$ is a normalized
angular measure on $r\T$.

The proof uses the following

\begin{claim}\label{claim3.2}
Given $\de \in (0, \frac14]$, $r\ge 1$, and
$z_0$, $\frac12 r \le |z_0| \le r$, there exists $\zeta\in z_0+
\de r\D$ such that
\[
\log|\psi(\zeta)| > \bigg(\frac12 - 3\de \bigg) |z_0|^2\,,
\]
unless $\omega$ belongs to an exceptional set of probability
$\exp\big( -c(\de)r^4\big)$.
\end{claim}

\begin{proof}[Proof of the claim]
The distribution (of probabilities) of the random potential
$\log|\psi(z)| - \frac12|z|^2$ is shift-invariant (see
\cite[Introduction]{SoTsi}). Writing the lower bound \eqref{eq2.3} in
Lemma~\ref{lem2.1} as
\[
\mathbb P \bigg( \max_{z\in r\D} \log|\psi(z)| - \tfrac12 |z|^2 \le
-\de r^2 \bigg) \le \exp\big( -c(\de) r^4 \big)
\]
we can apply it to the function $z \mapsto \log|\psi(z_0+z)| - \frac12
|z_0+z|^2$ on $ \de r \D $. We get 
\[
\mathbb P \bigg( \max_{z\in \de r\D} \log|\psi(z_0+z)| - \tfrac12
|z_0+z|^2 \le -\de (\de r)^2 \bigg) \le \exp\big( -c(\de) (\de r)^4
\big) \, .
\]
Assuming that $ \om $ does not belong to the exceptional set we obtain
$ z \in \de r \D $ such that
\[
\log|\psi(z+z_0)| - \tfrac12 |z+z_0|^2 \ge - \de^3 r^2 \, .
\]
Taking into account that $ |z| \le 2\de |z_0| $ we get $ \frac12
|z_0+z|^2 \ge \frac12 |z_0|^2 (1-2\de)^2 $;
\begin{multline*}
\log|\psi(z+z_0)| \ge \frac12 |z_0|^2 (1-2\de)^2 - \de^3 r^2 \\
\ge \frac12 |z_0|^2 - 2\de |z_0|^2 - \big( \tfrac14 \big)^2 \de
 (2|z_0|)^2 \ge \frac12 |z_0|^2 - 3\de |z_0|^2 \, ,
\end{multline*}
which yields the claim.
\end{proof}

\begin{proof}[Proof of Lemma \textup{\ref{lem3.1}}]
Now, we choose $\kappa = 1 - \de^{1/4}$, take
$N=[2\pi\de^{-1}]$, and consider $N$ discs (see Fig.~\ref{fig1})
\[
z_j+ \de r\D\,, \qquad z_j = \kappa r \exp\bigg( \frac{2\pi
ij}{N}\bigg)\,, \quad j=0, 1,\, ... \,, N-1\,.
\]
\begin{figure}
\begin{center}
\includegraphics{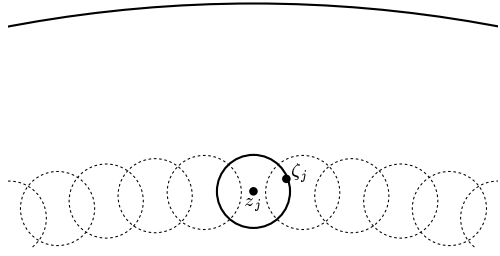}
\end{center}
\caption{\label{fig1}Small discs near the large circle}
\end{figure}
Claim~\ref{claim3.2} implies that if $\omega$ does not belong to
an exceptional set of probability $ N \exp(-c(\de)r^4) =
\exp(-c_1(\de) r^4)$, then we can choose $ N $ points $ \zeta_j
\in z_j + \de r \D$ such that
\[
\log|\psi(\zeta_j)| \ge \bigg( \frac12 - 3\de\bigg) |z_j|^2 \ge
\bigg( \frac12 - C \de^{1/4} \bigg) r^2\,.
\]
Let $P(z,\zeta)$ be the Poisson kernel for the disc $r\D$,
$|z|=r$, $|\zeta|<r$. We set $P_j(z) = P(z, \zeta_j)$.
Then
\begin{eqnarray*}
\bigg( \frac12 - C\de^{1/4} \bigg) r^2 &\le& \frac1{N}
\sum_{j=0}^{N-1} \log|\psi(\zeta_j)| \le
\int_{r\T} \bigg( \frac1{N} \sum_{j=0}^{N-1}P_j \bigg)
\log|\psi|\, d\mu \\
&=& \int_{r\T} \log|\psi|\, d\mu + \int_{r\T} \bigg(
\frac1{N}\sum_{j=0}^{N-1} P_j - 1\bigg) \log|\psi|\, d\mu\,.
\end{eqnarray*}
We have
\[
\int_{r\T} \bigg( \frac1{N}\sum_{j=0}^{N-1} P_j - 1\bigg)
\log|\psi|\, d\mu \le \max_{z\in r\T} \bigg|
\frac1{N}\sum_{j=0}^{N-1} P_j - 1\bigg|\, \cdot\, \int_{r\T}
\big| \log|\psi|\, \big|\, d\mu\,.
\]
The next two claims finish the job.
\end{proof}

\begin{claim}\label{claim3.3}
\[
\max_{z\in r\T} \bigg| \frac1{N}\sum_{j=0}^{N-1} P_j - 1\bigg| \le C
\de^{1/2}\,.
\]
\end{claim}
\begin{claim}\label{claim3.4}
\[
\int_{r\T} \big| \log|\psi| \,\big|\, d\mu \le 10 r^2
\]
provided that  $\omega$ does not belong to an exceptional set of
probability $\exp(-c r^4)$.
\end{claim}

\begin{proof}[Proof of Claim \textup{\ref{claim3.3}}]
We start with
\[
\int_{\kappa r\T} P(z, \zeta)\, d\mu(\zeta) = 1\,,
\]
and split the circle $\kappa r\T$ into a union of $N$ disjoint arcs
$I_j$ of equal angular measure $\mu(I_j)= \frac1{N}$ centered at
$z_j$. Then
\[
1 = \frac1{N} \sum_{j=0}^{N-1} P(z, \zeta_j)
+ \sum_{j=0}^{N-1} \int_{I_j} \big( P(z, \zeta)
- P(z, \zeta_j) \big) \, d\mu(\zeta)\,,
\]
and
\begin{eqnarray*}
|P(z, \zeta) - P(z, \zeta_j)| &\le& \max_{\zeta\in I_j} |\zeta-\zeta_j|
\, \cdot\, \max_{z, \zeta} |\nabla_\zeta P(z, \zeta)| \\
&\le& C_1\de r \, \cdot \, \frac{C_2 r}{(r-|\zeta|)^2} =
\frac{C\de}{\de^{1/2}} = C\de^{1/2}\,,
\end{eqnarray*}
proving the claim.
\end{proof}

\begin{proof}[Proof of Claim \textup{\ref{claim3.4}}]
By Lemma~\ref{lem2.1}, we know that unless $\omega$ belongs to an
exceptional set of probability $\exp(-c r^4)$, there is a point
$\zeta\in \frac12 r\T$ such that $\log |\psi(\zeta)| \ge 0$. Fix
such a $\zeta$. Then
\[
0 \le \int_{r\T} P(z, \zeta) \log|\psi(z)|\, d\mu(z)\,,
\]
and hence
\[
\int_{r\T} P(z, \zeta) \log^-|\psi(z)|\, d\mu(z) \le \int_{r\T}
P(z, \zeta) \log^+|\psi(z)|\, d\mu(z)\,.
\]
It remains to recall that for $|z|=r$ and $|\zeta|=\frac12 r$,
\[
\frac13 \le P(z, \zeta) \le 3\,,
\]
and that
\[
\int_{r\T} \log^+ |\psi|\, d\mu \le \log M(r, \psi) \le r^2
\]
(provided $\omega$ is non-exceptional). Hence
\[
\int_{r\T} \log^-|\psi|\, d\mu \le 9 r^2\,,
\]
and
\[
\int_{r\T} \big| \log|\psi| \,\big|\, d\mu \le 10 r^2\,,
\]
proving the claim.
\end{proof}

\section{Proof of Theorem \ref{thm2}}
We shall prove that
\begin{equation}\label{eq4.1}
\mathbb P \bigg( \frac{n(r)}{r^2} \le  1 + \de \bigg) \le
\exp\big( -c(\de)r^4\big)\,.
\end{equation}
The proof of the lower bound for $n(r)$ is practically the same
and is left to the reader.

\begin{sloppypar}
Fix $\kappa = 1 + \sqrt{\de}$. Then by Jensen's formula
\cite[Chapter~5, Section~3.1]{Ahlfors}
\[
n(r) \log\kappa \le \int_{r}^{\kappa r} \frac{n(t)}{t}\, dt
= \bigg( \int_{\kappa r\T} - \int_{r\T}\bigg) \log|\psi|\, d\mu\,,
\]
whence by Lemmas~\ref{lem2.1} and \ref{lem3.1}
\[
\frac{n(r)}{r^2} \le \frac1{\log\kappa}
\bigg( \kappa^2 \Big( \frac12 + \de \Big)
- \Big( \frac12 -\de \Big)\bigg) = \frac12 \frac{\kappa^2 -1}{\log\kappa}
+ \de \frac{\kappa^2+1}{\log\kappa} \le 1 + C\sqrt{\de}\,,
\]
provided that $\omega$ does not belong to an exceptional set of
probability $ {\exp \big( -c(\de) r^4\big)} $. This proves estimate
\eqref{eq4.1}.
\end{sloppypar}

\bigskip
\filbreak { \small
\begin{sc}
\parindent=0pt\baselineskip=12pt
\parbox{2.3in}{
Mikhail Sodin\\
School of Mathematics\\
Tel Aviv University\\
Tel Aviv 69978, Israel
\smallskip
\emailwww{sodin@tau.ac.il} {} } \hfill
\parbox{2.3in}{
Boris Tsirelson\\
School of Mathematics\\
Tel Aviv University\\
Tel Aviv 69978, Israel
\smallskip
\emailwww{tsirel@tau.ac.il} {www.tau.ac.il/\textasciitilde
tsirel/} }
\end{sc}
} \filbreak

\end{document}

%% file: preamble.tex
\numberwithin{equation}{section}

\newtheorem{theorem}{Theorem}
\newtheorem{lemma}{Lemma}
\newtheorem{claim}{Claim}

\newcommand{\const}{\operatorname{const}}

\newcommand{\om}{\omega}
\newcommand{\de}{\delta}

\newcommand{\R}{\mathbb R}
\newcommand{\C}{\mathbb C}
\newcommand{\T}{\mathbb T}
\newcommand{\D}{\mathbb D}

\def\emailwww#1#2{\par\qquad {\tt #1}\par\qquad {\tt #2}\medskip}